\documentclass[12pt, a4paper]{amsart} 

\usepackage{mathptmx} 
\usepackage{amsfonts,amsthm,amssymb,dsfont,stmaryrd,mathrsfs} 
\usepackage[leqno]{amsmath} 
\usepackage{array}
\usepackage{colortbl}
\usepackage{comment} 

\usepackage{graphicx,caption} 
\usepackage[pdfstartview=FitH]{hyperref} 
\usepackage{pdfpages} 
\usepackage{longtable,supertabular, multirow} 
\usepackage{doi}
\usepackage{tikz, tikz-cd}
\usetikzlibrary{matrix,arrows,decorations.pathmorphing}

\usepackage[hmarginratio={1:1}, vmarginratio={1:1}, top=3cm, bottom=2cm, left=3cm, right=3cm]{geometry} 

\usepackage{setspace}
\usepackage{mathtools} 

\theoremstyle{plain}
\newtheorem*{thm*}{Theorem} 
\newtheorem{thm}{Theorem} 
\newtheorem*{thm2*}{Theorem \ref{thmDecNilp}}
\newtheorem*{thm3*}{Theorem \ref{thmDecNilp2}}

\newtheorem{lem}[thm]{Lemma}
\newtheorem*{lem*}{Lemma}

\newtheorem*{prop*}{Proposition}

\newtheorem*{cor*}{Corollary}

\theoremstyle{definition}

\newtheorem*{defn*}{Definition}

\newtheorem*{conjecture*}{Conjecture}
\newtheorem{exmp}[thm]{Example}
\newtheorem*{exmp*}{Example}

\newtheorem*{prob*}{Problem}

\newtheorem*{ques*}{Question}

\newtheorem*{remk*}{Remark}

\def\nn{\mathbb{N}}

\def\qq{\mathbb{Q}}

\def\ff{\mathbb{F}}

\def\disp{\displaystyle}

\def\ord{{\rm ord}}

\begin{document}

\title{A note on prime index of a certain subgroup in $\ff_p^*$}

\author{Wei-Liang Sun}
\address{Department of Mathematics, National Taiwan Normal University, Taipei 11677, Taiwan, ROC}
\email{wlsun@ntnu.edu.tw}

\begin{abstract}
Under the generalized Riemann hypothesis, we illustrate that the ratio of the set of primes $p$ such that $\langle -1, 2 \rangle$ has an odd prime index in $\ff_p^*$ to the set of primes $p$ such that the subgroup has index greater than $2$ nears $46 \%$. 
\end{abstract}

\maketitle

Let $p$ be an odd prime number and let $\ff_p$ be the finite field consisting of $p$ elements. Set $\ell_p = [\ff_p^* : \langle -1, 2 \rangle]$, the group index of the subgroup generated by $-1$ and $2$ in the multiplicative group $\ff_p^* = \ff_p \setminus \{0\}$. For $\ell \in \nn$, let $d_{\ell}$ be the natural density of the set $\{\text{prime } p \mid \ell_p = \ell\}$. Consider the following ratio $$r = \frac{\disp \sum_{\ell: \text{ prime} \geq 3} d_{\ell}}{\disp \sum_{\ell \geq 3} d_{\ell}}.$$ We will obtain that $0.462 < r < 0.463$ under the generalized Riemann hypothesis. (By a computer search to the first hundred million primes, we see $r \approx 0.462535 \ldots$.) 
The consideration of this ratio $r$ is motivated by the article [Hsia-Li-Sun]\footnote{L.-C. Hsia, H.-C. Li and W.-L. Sun. {\it Conflict-avoiding codes of prime lengths and cyclotomic numbers}. \doi{10.48550/arXiv.2302.01487}.} 

First of all, we rewrite $\ell_p$ as $\frac{p-1}{2 \, \ord_p(4)}$ where $\ord_p(n)$ is the order of element $n$ in $\ff_p^*$. Indeed, if $\ord_p(2)$ is odd, then $| \langle -1, 2 \rangle| = 2 \ord_p(2) = 2 \ord_p(4)$; if $\ord_p(2)$ is even, then $| \langle -1, 2 \rangle| = \ord_p(2) = 2 \ord_p(4)$.\footnote{We thank Professor Pieter Moree for pointing out this fact to us.} Thus, $d_{\ell}$ is the natural density of the set of primes $p$ for which $4$ has residue index $2 \ell$ modulo $p$. It follows that $d_{\ell} = A(4, 2 \ell)$ in the terminology of [Wagstaff (1982)]\footnote{S.S. Wagstaff Jr., {\it Pseudoprimes and a generalization of Artin's conjecture}, Acta Arith. 41 (1982) 141-150. \doi{10.4064/aa-41-2-141-150}.}. Assume the generalized Riemann hypothesis for the number fields $\qq(\zeta_k, 4^{1/k})$ where $k$ is a positive integer and $\zeta_k$ is a primitive $k$-th root of unity. Then Example~3 of [Wagstaff (1982)] presents $$d_{\ell} = \left\{ \begin{array}{rl} \frac{3}{2} \, A \, Q(\ell) & \text{if $\ell$ is odd,} \\ \frac{1}{3} \, A \, Q(\ell) & \text{if $\ell \equiv 2 \pmod{4}$,} \\ A \, Q(\ell) & \text{if $4 \mid \ell$,} \end{array} \right.$$ where $$A = \prod_{p: \text{ prime}} \left( 1 - \frac{1}{p (p-1)} \right) = 0.3739558136...$$ is the Artin constant and $$Q(\ell) = \frac{1}{\ell^2} \prod_{q \mid \ell} \frac{q^2-1}{q^2-q-1}$$ where $q$ is a prime divisor of $\ell$. In particular, $d_{\ell} > 0$ for every $\ell$.

Note that the sum of densities of all situations must be $1$. So $\sum_{\ell \geq 1} d_{\ell} = 1$. Thus, the denominator of $r$ is \begin{align} \sum_{\ell \geq 3} d_{\ell} = 1 - d_1 - d_2 = 1 - \frac{3}{2} A Q(1) - \frac{1}{3} A Q(2) = 1 - \frac{7}{4} A. \label{sum of densities} \end{align}

For the numerator of $r$, one has \begin{align} \sum_{\ell: \text{ prime} \geq 3} d_{\ell} = \frac{3}{2} A \sum_{\ell: \text{prime} \geq 3} Q(\ell) = \frac{3}{2} A \left( \sum_{\ell: \text{ prime} \geq 3} \left(Q(\ell) - \frac{1}{\ell^2}\right) + P(2) - \frac{1}{4} \right) \end{align} where $P(s) = \sum_{\text{prime } p} \frac{1}{p^s}$ is the prime zeta function. 

Since $Q(\ell) > \frac{1}{\ell^2}$, it follows that \begin{align} r = \frac{\disp \sum_{\ell: \text{ prime} \geq 3} d_{\ell}}{\disp \sum_{\ell \geq 3} d_{\ell}} > L(m) :=\frac{\frac{3}{2} A \left( \disp \sum_{\ell: \text{ prime} \geq 3}^{m} \left(Q(\ell) - \frac{1}{\ell^2}\right) + P(2) - \frac{1}{4} \right)}{1 - \frac{7}{4} A} \end{align} for every positive integer $m \geq 3$.  The lower bound $L(m)$ of $r$ is increasing when $m$ grows. As $P(2) = 0.452247420041065 \ldots$, we have the following table of several values of $L(m)$:
$$\begin{array}{|c|c|c|c|c|c|} \hline
m & 3 & 5 & 7 & 11 & 13 \\ \hline
L(m) & 0.436495\ldots & 0.453581\ldots & 0.459237\ldots & 0.460591\ldots & 0.461396\ldots \\ \hline
\end{array}$$
$$\begin{array}{|c|c|c|c|c|c|} \hline
m & 17 & 19 & 23 & 29 & 31 \\ \hline
L(m) & 0.461749\ldots & 0.461999\ldots & 0.462139\ldots & 0.462208\ldots & 0.462264\ldots \\ \hline
\end{array}$$


For an upper bound of $r$, we obverse that $Q(\ell) - \frac{1}{\ell^2} = \frac{1}{\ell (\ell^2 - \ell - 1)} \leq \frac{49}{41} \, \frac{1}{\ell^3}$ for $\ell \geq 7$. Therefore \begin{align} r \leq \frac{\frac{3}{2} A \left( \disp \frac{1}{15} + \frac{1}{95} + \frac{49}{41} \left(P(3) - \frac{1}{2^3} - \frac{1}{3^3} - \frac{1}{5^3} \right) + P(2) - \frac{1}{4} \right)}{1 - \frac{7}{4} A} = 0.462748 \ldots.\end{align} where $P(3) = 0.174762639299443 \ldots$.

\end{document}